\documentclass[11pt,reqno]{amsart}
\usepackage[a4paper,margin=28mm]{geometry}
\usepackage{amsmath,amssymb,amsfonts,amsthm,mathtools}
\usepackage{microtype}
\usepackage[colorlinks=true,linkcolor=blue,citecolor=blue,urlcolor=blue]{hyperref}
\numberwithin{equation}{section}
\newtheorem{theorem}{Theorem}[section]
\newtheorem{proposition}{Proposition}[section]

\newtheorem{remark}[theorem]{Remark}
\newcommand{\Ric}{\operatorname{Ric}}
\newcommand{\diver}{\operatorname{div}}
\newcommand{\tr}{\operatorname{tr}}
\newcommand{\II}{\mathrm{II}}
\title[Alexandrov's theorem in the anti-Gaussian space]{Alexandrov's theorem in the Anti-Gaussian Space: \\A Reilly-formula approach}

\author[Bai]{Jinchuan Bai}
\address[J.B.]{School of Mathematical Sciences\\
Xiamen University\\
361005\\
Xiamen\\
P. R.  China}
\email{baijinchuan@stu.xmu.edu.cn
}

\author[Xia]{Chao Xia}
\address[C.X.]{School of Mathematical Sciences\\
Xiamen University\\
361005\\
Xiamen\\
P. R.  China}
\email{chaoxia@xmu.edu.cn}
\thanks{This work is supported by
NSFC (Grant No.  12271449,  12671069, 12526203, 12526102) and the Natural Science Foundation of
Fujian Province of China (Grant No.  2024J011008).}

\date{}
\begin{document}
\begin{abstract}
We give an alternative proof of Alexandrov's theorem about constant weighted mean curvature hypersurfaces in the anti-Gaussian space, due to the authors \cite{BaiXia} and Johne-Silini \cite{JohneSilini} independently, via a 
Reilly-type formula approach. Moreover, we establish a sharp Heintze-Karcher-type inequality for hypersurfaces in manifolds with a generalized substatic condition, together with its rigidity in Euclidean space with density.
\end{abstract}

\maketitle

\section{Introduction}\label{sec:intro}

Alexandrov's soap bubble theorem asserts that every closed, connected, embedded
hypersurface of constant mean curvature in Euclidean space is a
round sphere \cite{Alexandrov}. Alexandrov's original proof uses his moving plane method. An alternative proof, which is based on the Heintze--Karcher inequality and the Minkowski formula, was later given by Reilly \cite{Reilly} and Ros \cite{Ros}.
The Heintze--Karcher inequality has been first established by Heintze-Karcher using
the normal exponential map \cite{HeintzeKarcher}, while
Reilly's integral formula provides another approach through a boundary value problem \cite{Reilly, Ros}. 

Heintze--Karcher's inequality and Alexandrov's theorem has been studied in other space forms \cite{MontielRos, Korevaar}. Moreover, Alexandrov's theorem for a class of substatic warped product manifolds has been proved by Brendle \cite{Brendle}, via establishment of Heintze--Karcher's inequality by means of the normal exponential map with respect to a conformal metric. 
We recall the concept of sub-staticity.  A Riemannian triple $(M, \overline g, V)$ is called sub-static if $V\in C^\infty(M)$ such that $V\ge 0$ in $M$ with $V=0$ only on $\partial M$, and
\begin{equation}
\overline{\Delta} V \overline g-\overline{\nabla}^2 V +V\overline{\rm Ric}\ge 0.
\end{equation}
Li-Xia \cite{LiXia} (see also Qiu-Xia \cite{QiuXia}) has established a generalized Reilly-type formula which fits perfectly for substatic manifolds and used it to generalize Brendle's Heintze-Karcher inequality to general substatic manifolds. Later, Borghini-Fogagnolo-Pinamonti \cite{BorghiniFogagnoloPinamonti} further proved the rigidity for such inequality.

Manifolds with density attracts many attention in recent decades. Let $(M^{m+1},\bar g)$ be a
Riemannian manifold endowed with a smooth weighted measure
$
d\mu_\phi=e^{-\phi}\,d\mu_{\bar g},
$
where $\phi\in C^\infty(M)$ is called the potential function. If
$\Omega\subset M$ is a domain with boundary
$\Sigma=\partial\Omega$, its weighted volume and the weighted area of
$\Sigma$ are defined, respectively, by
\[
\operatorname{Vol}_\phi(\Omega)
=
\int_\Omega d\mu_{\phi},
\qquad
\operatorname{Area}_\phi(\Sigma)
=\int_\Sigma d\sigma_{\phi}=
\int_\Sigma e^{-\phi}\,d\sigma.
\]
A particularly important model is Euclidean space equipped with a radial log-convex density. Manifolds
with density and their associated isoperimetric problems have been
extensively investigated; see, for example, \cite{Morgan2005,
RosalesCaneteBayleMorgan2008}.
In particular, the Log-Convex Density Conjecture, commonly attributed to Brakke (see \cite{
RosalesCaneteBayleMorgan2008}), has been resolved in full by Chambers
\cite{Chambers2019}, which asserts that balls centered at the origin are the only isoperimetric regions in $\mathbb{R}^{m+1}$ with radial log-convex density.

The stationary hypersurfaces for the isoperimetric problems in such spaces are given by constant weighted mean curvature hypersufaces $$H_\phi=H-\langle\overline \nabla \phi,\nu\rangle=\lambda,$$ where $H$ is the sum of the principal curvatures, $\nu$ is
the outward unit normal,  $\lambda\in \mathbb{R}$.
In an important special case $\phi=-|x|^2/4$, the constant weighted mean curvature hypersufaces satisify
\[
 H_\phi=H+\frac12\langle x,\nu\rangle=\lambda,
\]
It is also called the $\lambda$-self-expander
equation, compare to the $\lambda$-self-shrinkers as  stationary hypersurfaces in the Gaussian space.

 We \cite{BaiXia} and Johne-Silini \cite{JohneSilini} have independently proved the following Alexandrov-type theorem for $\lambda$-self-expanders.
 \begin{theorem}[Bai-Xia \cite{BaiXia}, Johne-Silini \cite{JohneSilini}]\label{thm:alex}
Every closed, connected, embedded $\lambda$-self-expander in $\mathbb{R}^{m+1}$ is a sphere
centered at the origin.
\end{theorem}
There are non-compact and immersed non-spherical $\lambda$-self-expander, \cite{AncariCheng2023,Qi2026}, and also several rigidity results of $\lambda$-self-expanders, see for example \cite{AncariCheng2023,BaiXiaStability,DengYaoYuZhang2024,LiQi2026}. 

The proof of Theorem \ref{thm:alex} is based on the following Minkowski formula for $\lambda$-self-expanders
\begin{equation}\label{eq:anti-Minkowski}
 \int_\Sigma \left(m+\frac{|x|^2}{2}\right)
 e^{|x|^2/4}\,d\sigma
= H_\phi \int_\Omega\left(m+1+\frac{|x|^2}{2}\right)
 e^{|x|^2/4}\,dx,
\end{equation}
and the following sharp Heintze-Karcher-type inequality.
 \begin{theorem}[Bai-Xia \cite{BaiXia}, Johne-Silini \cite{JohneSilini}]\label{thm:anti}
Let $\Omega\subset\mathbb R^{m+1}$ be a bounded domain with boundary $\Sigma=\partial\Omega$. If
$H+\langle x,\nu\rangle/2>0$ on $\Sigma$, then
\begin{equation}\label{eq:anti-HK}
 \int_\Sigma\frac{m+|x|^2/2}{H+\langle x,\nu\rangle/2}
 e^{|x|^2/4}\,d\sigma
 \ge\int_\Omega\left(m+1+\frac{|x|^2}{2}\right)
 e^{|x|^2/4}\,dx.
\end{equation}
Equality holds if and only if $\Omega$ is a ball centered at the origin.
\end{theorem}
The proof of Theorem \ref{thm:anti} is based on an adaption of Brendle's idea by using normal flow with respect to the conformal metric $(m+|x|^2/2)^{-2}\overline g$.
We remark that the approach also works in a large class of warped product manifolds with radial density satisfying a $\phi$-substatic assumption, see the details in \cite{BaiXia}.

In this paper, inspired by Li-Xia's work \cite{LiXia}, we find a Reilly-formula-approach to Theorem \ref{thm:anti}, and more generally, to Heintze-Karcher-type inequality in $\phi$-substatic manifolds (without warped product structure for metric and radial structure for density).
To state our result, we introduce some notations. Let
$\overline \Delta_\phi =\overline\Delta
       -\langle\overline\nabla\phi,\overline\nabla \rangle$ be the Witten-Laplacian and 
$\overline\Ric_\phi=\overline\Ric+\overline\nabla^2\phi$ be the Bakry-Emery-Ricci curvature. Denote
\begin{equation}\label{eq:weighted-ric-intro}
 \overline\Ric^{V}_{\phi}
 :=\frac{\overline \Delta_\phi V}{V}\overline g
   -\frac{\overline\nabla^2V}{V}+\overline\Ric_{\phi},
\end{equation}
and we call it $\phi$-static-Ricci tensor. Our main result is the following.
\begin{theorem}\label{thm:general}
Let $(M^{m+1},\overline g,e^{-\phi}d\mu_{\overline g})$ be a Riemannian manifold with density, and  $\Omega\subset M$
be a bounded domain with boundary
$\Sigma=\partial\Omega$ such that $H_\phi>0$. Let $V\in C^\infty(\overline\Omega)$ and $\mathcal N\in C^\infty(\overline\Omega)$ be such that $V>0$, $\mathcal N\ge m+1$ and
\begin{equation}\label{eq:tensor-condition}
 (\mathcal N-(m+1))\overline\Ric^{V}_{\phi}
 \ge\overline\nabla\phi\otimes\overline\nabla\phi.
\end{equation}
At points where $\mathcal N=m+1$, assume in addition that $\overline\nabla\phi=0$ and $\overline\Ric^V_\phi\geq0$. 
Then it holds that
\begin{equation}\label{eq:general-HK}
 \int_\Sigma\frac{V}{H_\phi}\,d\sigma_\phi
 \ge\int_\Omega\frac{\mathcal N}{\mathcal N-1}
 V\,d\mu_\phi.
\end{equation}
\end{theorem}
\begin{remark}
When $\phi\equiv 0$, we can choose $\mathcal{N}\equiv m+1$ and Theorem \ref{thm:general} reduces to Li-Xia's \cite{LiXia} Heintze-Karcher inequality for substatic manifolds. For general $\phi$ and a constant $\mathcal{N}>m+1$, Theorem \ref{thm:general} reduces to Huang-Ma-Zhu \cite{HMZ}. 
\end{remark}

The proof of Theorem \ref{thm:general} uses Li-Xia's  Reilly-type formula \cite{LiXia}, see also Huang-Ma-Zhu \cite{HMZ}. For simplicity, we state here a special case when $f|_{\Sigma}=0$.
\setcounter{proposition}{4}
\begin{proposition}[Li-Xia \cite{LiXia}, Huang-Ma-Zhu \cite{HMZ}]\label{prop:reilly}
Let $f\in C^\infty(\overline\Omega)$ be such that $f|_{\Sigma}=0$.
 Then
\begin{equation}\label{eq:reilly-intro}
\begin{aligned}
 \int_\Sigma VH_\phi f_\nu^2\,d\sigma_\phi
 ={}&\int_\Omega V\biggl[
 \left(\overline \Delta_\phi f-\frac fV\overline \Delta_\phi V\right)^2
 -\left|\overline\nabla^2 f-\frac fV\overline\nabla^2V\right|^2\biggr]
 \,d\mu_\phi\\
 &-\int_\Omega V\overline\Ric^V_\phi
 \left(\overline\nabla f-\frac fV\overline\nabla V,
       \overline\nabla f-\frac fV\overline\nabla V\right)
 \,d\mu_\phi.
\end{aligned}
\end{equation}
\end{proposition}

Next we consider the rigidity for the Heintze-Karcher inequality \eqref{eq:general-HK} in Theorem \ref{thm:general}. For this purpose, we need a slightly stronger condition \eqref{eq:tensor-condition-strong} than \eqref{eq:tensor-condition}.
For simplicity, we only consider the case $M=\mathbb{R}^{m+1}$ and choose
\begin{equation}\label{NV}
 \mathcal{N}=m+1-\langle x,\overline\nabla\phi\rangle,
\qquad V=\mathcal{N}-1=m-\langle x,\overline\nabla\phi\rangle.
\end{equation}

To ensure $\mathcal N\ge m+1$, we assume
$\langle x,\overline\nabla\phi\rangle\le0$.



\setcounter{theorem}{5}
\begin{theorem}\label{cor:intro-radial}
Let $\Omega\subset\mathbb R^{m+1}$ be a bounded
connected smooth domain with density $d\mu_\phi=e^{-\phi}dx$
and boundary $\Sigma=\partial\Omega$ such that $H_\phi>0$.
 Choose $\mathcal N$ and $V$ as in \eqref{NV}. Assume that $\langle x,\overline\nabla\phi\rangle\le 0$, 
and \begin{equation}\label{eq:tensor-condition-strong}
\overline\Ric^V_\phi\ge
\frac{|\overline\nabla\phi|^2}{\mathcal N-(m+1)}\overline g,
\end{equation}
where the quotient on the right hand side admits a smooth, strictly positive
extension to $\overline\Omega$.
Then \eqref{eq:general-HK} holds. Moreover, equality in \eqref{eq:general-HK} holds if and only if $\Omega=B_R(0)$ for some $R>0$ and 
$\phi$ is a radial function about the origin.
\end{theorem}

\begin{remark}
Suppose that $\phi=\phi(r)$, where $r=|x|$.
Since
\[
\frac{|\overline\nabla\phi|^2}
{\mathcal N-(m+1)}
=-\frac{\phi'(r)}{r},
\]
condition \eqref{eq:tensor-condition-strong}
reduces to
\[
\overline\Delta_\phi V\,\overline g
-\overline\nabla^2V
+V\overline\Ric_\phi
+V\frac{\phi'(r)}{r}\overline g
\geq 0.
\]
This is precisely the weighted substatic
condition used in \cite{BaiXia} for the
Euclidean case $\lambda(r)=r$. In this case, together with the smooth positive
extension assumption, this condition is stronger
than radial log-convexity on a ball centered at
the origin(cf. \cite[Appendix B]{JohneSilini}).
\end{remark}

We return to the case of the anti-Gaussian space, that is, $\phi=-\frac{|x|^2}{4}$. Then 
\begin{equation*}V=m+\frac{|x|^2}{2},\qquad 
\mathcal N=m+1+\frac{|x|^2}{2}\end{equation*}
and
\begin{equation*}
(\mathcal N-(m+1))\overline\Ric^{V}_{\phi}=\frac{|x|^2}{4}\overline g= |\overline\nabla\phi|^2\overline g\ge\overline\nabla\phi\otimes\overline\nabla\phi.
\end{equation*}

Moreover,
\[
\overline\Ric^V_\phi
=\frac12\overline g
=\frac{|\overline\nabla\phi|^2}
{\mathcal N-(m+1)}\overline g,
\]
where the quotient extends smoothly to the origin.

Hence, Theorem \ref{thm:anti} can be viewed as a special case of Theorem \ref{thm:general} and Theorem \ref{cor:intro-radial}.

Finally, we mention a related work by de Lima \cite{delima}, where an Alexandrov-type theorem in the Gaussian half-space $(\mathbb{R}^{m+1}_+, e^{-\frac{|x|^2}{4}}dx)$ has been proved via Qiu-Xia's Reilly-type formula \cite{QiuXia} with $\infty$-Bakry-Emery-Ricci. Compare to de Lima \cite{delima}, the Alexandrov-type theorem in the anti-Gaussian space proved here uses Li-Xia's Reilly-type formula \cite{LiXia} with $\phi$-static-Ricci and a special function $\mathcal{N}$.
    
The rest of the paper is organized as follows. In Section~\ref{sec:formula} we give a proof of Proposition \ref{prop:reilly} for convenience of reader and then prove Theorem \ref{thm:general}.
In Section~\ref{sec:positive}, we  prove Theorem \ref{cor:intro-radial}.

\subsection*{Acknowledgments}
Li–Xia's Reilly formula used in this paper was obtained while CX was a postdoctoral fellow at McGill University under the mentorship of Professor Pengfei Guan. CX is deeply grateful to Professor Pengfei Guan for his continued support.

\section{A weighted Reilly formula and a Heintze--Karcher inequality}\label{sec:formula}

Throughout this section we use the notation of the introduction.
The unbarred operators $\nabla$, $\Delta$, and $\diver$ act on
$\Sigma$, and for ambient vector fields $X$ we write
\[
 \overline\diver_\phi X
 =\overline\diver X-\langle\overline\nabla\phi,X\rangle.
\]
Our convention is $\II(X,Z)=\langle\overline\nabla_X\nu,Z\rangle$
for $X,Z\in T\Sigma$, with $H=\tr_\Sigma\II$.

Related Reilly-type integral formulas were established by
Qiu--Xia \cite{QiuXia} and Li--Xia \cite{LiXia}, and extended
to manifolds with density by Huang--Ma--Zhu \cite{HMZ}.
For completeness, we give a brief proof of
Proposition~\ref{prop:reilly}.
\begin{proof}[Proof of Proposition~\ref{prop:reilly}]
Set
\[
 Y=\overline\nabla f-\frac fV\overline\nabla V
   =V\overline\nabla\left(\frac fV\right).
\]
Direct differentiation gives
\begin{equation}\label{eq:reilly-Y-derivatives}
\begin{aligned}
 \overline\diver_\phi(VY)
 &=V\overline \Delta_\phi f-f\overline \Delta_\phi V,\\
 \overline\nabla_XY
 &=\overline\nabla_X\overline\nabla f
   -\frac fV\overline\nabla_X\overline\nabla V
   -\frac{\langle X,Y\rangle}{V}\overline\nabla V.
\end{aligned}
\end{equation}
It follows that
\begin{align}
&\overline\diver_\phi
 \bigl((V\overline \Delta_\phi f-f\overline \Delta_\phi V)Y\bigr)
 \notag\\
={}&V\left(\overline \Delta_\phi f-\frac fV\overline \Delta_\phi V\right)^2
 +V\left\langle
 \overline\nabla\left(\overline \Delta_\phi f
 -\frac fV\overline \Delta_\phi V\right),Y\right\rangle.
\label{eq:reilly-div-first}
\end{align}
To calculate the second divergence, use the weighted
Bochner commutation formula
\[
 \overline\diver_\phi(\overline\nabla^2h)(Z)
 =\langle\overline\nabla\overline \Delta_\phi h,Z\rangle
  +\overline\Ric_\phi(\overline\nabla h,Z),
\]
where the divergence of a Hessian is taken in its first
index. Applied to $f$ and $V$, this formula and
\eqref{eq:reilly-Y-derivatives} give
\begin{align*}
&\overline\diver_\phi
 \bigl(V\overline\nabla_Y\overline\nabla f
  -f\overline\nabla_Y\overline\nabla V\bigr)\\
={}&V\langle\overline\nabla\overline \Delta_\phi f,Y\rangle
 -f\langle\overline\nabla\overline \Delta_\phi V,Y\rangle
 +V\overline\Ric_\phi(Y,Y)\\
&+V\left|\overline\nabla^2f
   -\frac fV\overline\nabla^2V\right|^2
 -\overline\nabla^2V(Y,Y).
\end{align*}
Since
\[
\begin{aligned}
&V\langle\overline\nabla\overline \Delta_\phi f,Y\rangle
 -f\langle\overline\nabla\overline \Delta_\phi V,Y\rangle\\
={}&V\left\langle\overline\nabla\left(
 \overline \Delta_\phi f-\frac fV\overline \Delta_\phi V\right),Y\right\rangle
 +(\overline \Delta_\phi V)|Y|^2,
\end{aligned}
\]
the definition of $\overline\Ric_\phi^V$ reduces this to
\begin{align}
&\overline\diver_\phi
 \bigl(V\overline\nabla_Y\overline\nabla f
   -f\overline\nabla_Y\overline\nabla V\bigr)
 \notag\\
={}&V\left\langle\overline\nabla\left(
 \overline \Delta_\phi f-\frac fV\overline \Delta_\phi V\right),Y\right\rangle
 +V\left|\overline\nabla^2f-\frac fV\overline\nabla^2V\right|^2
 +V\overline\Ric_\phi^V(Y,Y).
\label{eq:reilly-div-second}
\end{align}

Subtracting \eqref{eq:reilly-div-second} from
\eqref{eq:reilly-div-first} and integrating gives the
interior expression in \eqref{eq:reilly-intro}. Its
boundary term is
\[
\begin{aligned}
\big\langle
 &(V\overline \Delta_\phi f-f\overline \Delta_\phi V)Y
 -V\overline\nabla_Y\overline\nabla f
 +f\overline\nabla_Y\overline\nabla V,\nu\big\rangle.
\end{aligned}
\]
Because $f=0$ on $\Sigma$, we have
$Y=\overline\nabla f=f_\nu\nu$. For tangent vectors
$X,Z\in T\Sigma$, differentiating $Z(f)=0$ gives
\[
 \overline\nabla^2f(X,Z)
 =-f_\nu\langle\nu,\overline\nabla_XZ\rangle
 =f_\nu\II(X,Z).
\]
Taking the tangential trace, and using the definition
of $\overline \Delta_\phi$, yields
\[
 \overline \Delta_\phi f-\overline\nabla^2f(\nu,\nu)
 =\bigl(H-\langle\overline\nabla\phi,\nu\rangle\bigr)f_\nu
 =H_\phi f_\nu.
\]
Thus the boundary integrand is $VH_\phi f_\nu^2$,
and the weighted divergence theorem proves
\eqref{eq:reilly-intro}.
\end{proof}

\begin{proof}[Proof of Theorem~\ref{thm:general}]
Let $f$ be the solution to the Dirichlet problem: 
\begin{equation}\label{eq:dirichlet-equation}
 \begin{cases}
 \displaystyle \overline{\Delta}_\phi f-\frac{\overline{\Delta}_\phi V}{V}f
 =\frac{\mathcal N}{\mathcal N-1}&\text{in }\Omega,\\[3pt]
 f=0&\text{on }\Sigma.
 \end{cases}
\end{equation}
Writing $f=Vu$, we have
\[
 \overline{\Delta}_\phi(Vu)-u\overline{\Delta}_\phi V
 =\frac1V\overline{\diver}_\phi(V^2\overline{\nabla} u),
\]
so the equation can be written in divergence form:
\[
 \overline{\diver}_\phi(V^2\overline{\nabla} u)
 =\frac{\mathcal N}{\mathcal N-1}V,
 \qquad u|_\Sigma=0.
\]
By the weighted divergence
theorem, we have
\begin{equation}\label{eq:boundary-integral}
 \int_\Sigma Vf_\nu\,d\sigma_\phi
 =\int_\Omega\frac{\mathcal N}{\mathcal N-1}V\,d\mu_\phi.
\end{equation}

The vector field from Proposition~\ref{prop:reilly} is
\[
 Y=\overline{\nabla}f-\frac fV\overline{\nabla}V.
\]
The trace estimate for the Hessian difference is
\[
\left|
\overline{\nabla}^{2}f
-\frac{f}{V}\overline{\nabla}^{2}V
\right|^{2}
\geq
\frac{1}{m+1}
\left(
\overline{\Delta}f
-\frac{f}{V}\overline{\Delta}V
\right)^{2}.
\]
At points where $\mathcal N>m+1$, the second tensor inequality
in \eqref{eq:tensor-condition} gives
\[
\overline{\Ric}^{V}_{\phi}(Y,Y)
\geq
\frac{1}{\mathcal N-(m+1)}
\langle\overline{\nabla}\phi,Y\rangle^{2}.
\]
Combining the two estimates with the scalar Cauchy--Schwarz
inequality yields
\begin{align}
&\left|
\overline{\nabla}^{2}f
-\frac{f}{V}\overline{\nabla}^{2}V
\right|^{2}
+\overline{\Ric}^{V}_{\phi}(Y,Y)
\nonumber\\
\geq&
\frac{1}{m+1}
\left(
\overline{\Delta}f
-\frac{f}{V}\overline{\Delta}V
\right)^{2}
+
\frac{1}{\mathcal N-(m+1)}
\langle\overline{\nabla}\phi,Y\rangle^{2}
\nonumber\\
\geq&
\frac{1}{\mathcal N}
\left(
\overline{\Delta}f
-\frac{f}{V}\overline{\Delta}V
-\langle\overline{\nabla}\phi,Y\rangle
\right)^{2}
\nonumber\\
=&
\frac{1}{\mathcal N}
\left(
\overline\Delta_{\phi}f
-\frac{f}{V}\overline\Delta_{\phi}V
\right)^{2}.
\label{eq:pointwise-estimate}
\end{align}
Since $\overline{\nabla}\phi=0$ and $\overline{\Ric}^{V}_{\phi}\geq0$ at points where $\mathcal N=m+1$,  
we also have \eqref{eq:pointwise-estimate}. Inserting the pointwise bound into the
Reilly identity of Proposition~\ref{prop:reilly}, and using
\eqref{eq:dirichlet-equation}, yields
\begin{align}
 \int_\Sigma VH_\phi f_\nu^2\,d\sigma_\phi
 &\le\int_\Omega V\left(1-\frac1{\mathcal N}\right)
       \left(\overline{\Delta}_\phi f-\frac{\overline{\Delta}_\phi V}{V}f\right)^2 d\mu_\phi\notag\\
 &=\int_\Omega\frac{\mathcal N}{\mathcal N-1}V\,d\mu_\phi.
 \label{eq:interior-inequality}
\end{align}
On $\Sigma$, \eqref{eq:boundary-integral} and Cauchy--Schwarz give
\begin{equation}\label{eq:boundary-cauchy}
\begin{aligned}
 \left(\int_\Omega\frac{\mathcal N}{\mathcal N-1}V\,d\mu_\phi\right)^2
 &=\left(\int_\Sigma Vf_\nu\,d\sigma_\phi\right)^2\\
 &\le\left(\int_\Sigma\frac V{H_\phi}\,d\sigma_\phi\right)
       \left(\int_\Sigma VH_\phi f_\nu^2\,d\sigma_\phi\right).
\end{aligned}
\end{equation}
Therefore,  \eqref{eq:interior-inequality} and \eqref{eq:boundary-cauchy} yields \eqref{eq:general-HK}.
\end{proof}

\section{Rigidity for Heintze-Karcher inequality}\label{sec:positive}

\begin{proof}[Proof of Theorem~\ref{cor:intro-radial}]
Since $\mathcal N=V+1$, the curvature condition
\eqref{eq:tensor-condition-strong} and the smooth positive
extension of $|\overline\nabla\phi|^2/(V-m)$ ensure that
all the assumptions of Theorem~\ref{thm:general} are satisfied,
including at $V=m$. Thus \eqref{eq:general-HK} follows.

Suppose equality holds, and take the solution $f$ of
\eqref{eq:dirichlet-equation}. As before, write
\[
u=\frac fV,\qquad
Y=\overline\nabla f-u\overline\nabla V
=V\overline\nabla u.
\]
Since $\mathcal N=V+1$, the Dirichlet equation gives
\[
\overline\Delta_\phi f-u\overline\Delta_\phi V
=\frac{V+1}{V}.
\]
By the definition of the weighted Laplacian,
\begin{equation}\label{eq:radial-trace}
\overline\Delta f-u\overline\Delta V
=\frac{V+1}{V}+\langle\overline\nabla\phi,Y\rangle.
\end{equation}

At points where $V>m$, the trace inequality and
\eqref{eq:tensor-condition-strong} yield
\begin{align*}
&\left|\overline\nabla^2f-u\overline\nabla^2V\right|^2
 +\overline\Ric_\phi^V(Y,Y)\\
\ge&
 \frac{(\overline\Delta f-u\overline\Delta V)^2}{m+1}
 +\frac{|\overline\nabla\phi|^2}{V-m}|Y|^2\\
\ge&
 \frac{(\overline\Delta f-u\overline\Delta V)^2}{m+1}
 +\frac{\langle\overline\nabla\phi,Y\rangle^2}{V-m}\\
\ge& \frac1{V+1}\bigl(
 \overline\Delta f-u\overline\Delta V
 -\langle\overline\nabla\phi,Y\rangle\bigr)^2\\
=&\frac{V+1}{V^2}.
\end{align*}
The second and third inequalities follow from the vector
and scalar Cauchy--Schwarz inequalities, respectively.
As in the proof of Theorem~\ref{thm:general}, equality
in \eqref{eq:general-HK} forces equality in these
pointwise estimates and in the boundary Cauchy--Schwarz
inequality.

Equality in the scalar Cauchy--Schwarz inequality gives
\[
\frac{\overline\Delta f-u\overline\Delta V}{m+1}
=-\frac{\langle\overline\nabla\phi,Y\rangle}{V-m}.
\]
Together with \eqref{eq:radial-trace}, this yields
\[
\begin{aligned}
\frac{V+1}{V}
&=\overline\Delta f-u\overline\Delta V
 -\langle\overline\nabla\phi,Y\rangle\\
&=\left(1+\frac{V-m}{m+1}\right)
 (\overline\Delta f-u\overline\Delta V)\\
&=\frac{V+1}{m+1}
 (\overline\Delta f-u\overline\Delta V).
\end{aligned}
\]
Consequently,
\begin{equation}\label{eq:radial-trace-equality}
\overline\Delta f-u\overline\Delta V=\frac{m+1}{V},
\qquad
\langle\overline\nabla\phi,Y\rangle=-\frac{V-m}{V}.
\end{equation}
Equality in the trace inequality and
\eqref{eq:radial-trace-equality} give
\[
\overline\nabla^2f-u\overline\nabla^2V=\frac1V\overline g.
\]
Equality in the vector Cauchy--Schwarz inequality
makes $\overline\nabla\phi$ and $Y$ linearly dependent.
Using the second identity in
\eqref{eq:radial-trace-equality}, we obtain
\[
\overline\nabla\phi
=-\frac{|\overline\nabla\phi|^2}{V-m}\,VY.
\]

At points where $V=m$, the smooth positive extension
of the quotient gives $\overline\nabla\phi=0$, and
\eqref{eq:tensor-condition-strong} implies
$\overline\Ric_\phi^V>0$. In this case
\eqref{eq:radial-trace} gives
$\overline\Delta f-u\overline\Delta V=(m+1)/m$, so
\[
\begin{aligned}
&\left|\overline\nabla^2f-u\overline\nabla^2V\right|^2
 +\overline\Ric_\phi^V(Y,Y)\\
&\qquad\ge\frac{m+1}{m^2}+\overline\Ric_\phi^V(Y,Y)
\ge\frac{m+1}{m^2}.
\end{aligned}
\]
Equality forces $Y=0$ and equality in the trace
inequality, hence the two interior identities remain
valid at these points as well, with the quotient
understood as its smooth extension.

Finally, equality in the boundary Cauchy--Schwarz
inequality implies that $H_\phi f_\nu$ is constant
on $\Sigma$. By \eqref{eq:boundary-cauchy} and equality
in \eqref{eq:general-HK}, this constant is $1$.
We have therefore obtained
\begin{equation}\label{eq:radial-equality}
\begin{gathered}
\overline\nabla^2f-u\overline\nabla^2V
=\frac1V\overline g,\\
\overline\nabla\phi
=-\frac{|\overline\nabla\phi|^2}{V-m}\,VY
\qquad\text{in }\Omega,\\
f_\nu=\frac1{H_\phi}
\qquad\text{on }\Sigma.
\end{gathered}
\end{equation}

We now prove that $\Omega$ is a ball centered at the
origin. Differentiating
$Y=\overline\nabla f-u\overline\nabla V$ and using
the first identity in \eqref{eq:radial-equality}, we get
\[
\overline\nabla_XY=\frac1V X-X(u)\overline\nabla V
\]
for every vector field $X$. Since
$Y=V\overline\nabla u$ and $\overline\nabla_Xx=X$,
\begin{align*}
\overline\nabla_X(VY-x)
&=X(V)Y+V\overline\nabla_XY-X\\
&=\langle\overline\nabla V,X\rangle Y
 -\langle Y,X\rangle\overline\nabla V.
\end{align*}
For any vector fields $X,Z$, the preceding identity gives
\[
\langle\overline\nabla_X(VY-x),Z\rangle
+\langle\overline\nabla_Z(VY-x),X\rangle=0.
\]
Thus $VY-x$ is a Euclidean Killing vector field.

By \eqref{eq:dirichlet-equation} and $\mathcal N=V+1$,
\[
\overline\diver_\phi(V^2\overline\nabla u)
=V+1>0,\qquad u|_\Sigma=0.
\]
The maximum principle gives $u<0$ in $\Omega$.
Thus $u$ attains an interior minimum at some
$x_0\in\Omega$, where $Y(x_0)=V\overline\nabla u(x_0)=0$.
The preceding formula then implies
$\overline\nabla_X(VY-x)(x_0)=0$ for every $X$.
Every Euclidean Killing vector field has the form
$Ax+b$, where $A$ is a constant skew-symmetric matrix
and $b$ is a constant vector. Since its derivative
vanishes at $x_0$, we have $A=0$. Evaluating at $x_0$,
where $Y(x_0)=0$, gives $b=-x_0$. Consequently,
\begin{equation}\label{eq:radial-center}
VY=x-x_0\qquad\text{in }\Omega.
\end{equation}

Substituting \eqref{eq:radial-center} into the second
identity in \eqref{eq:radial-equality}, we find
\[
\overline\nabla\phi
=-\frac{|\overline\nabla\phi|^2}{V-m}(x-x_0).
\]
In particular, $\overline\nabla\phi(x_0)=0$,
so \eqref{NV} gives $V(x_0)=m$.
As $V\ge m$, the function $V$ attains its minimum
at $x_0$ and hence $\overline\nabla V(x_0)=0$.
By \eqref{NV},
\[
V-m=\frac{|\overline\nabla\phi|^2}{V-m}
\langle x,x-x_0\rangle.
\]
Since
$\overline\nabla\langle x,x-x_0\rangle=2x-x_0$,
differentiating gives
\[
\begin{aligned}
\overline\nabla V={}
\langle x,x-x_0\rangle\,
\overline\nabla\left(
\frac{|\overline\nabla\phi|^2}{V-m}\right)+\frac{|\overline\nabla\phi|^2}{V-m}(2x-x_0).
\end{aligned}
\]
At $x=x_0$, the first term vanishes, so
\[
0=\overline\nabla V(x_0)
=\left.\frac{|\overline\nabla\phi|^2}{V-m}
\right|_{x_0}x_0.
\]
The quotient has a strictly positive smooth extension,
which forces $x_0=0$. Hence
\begin{equation}\label{eq:radial-position}
VY=x\qquad\text{in }\Omega.
\end{equation}

On $\Sigma$, we have $u=0$, and the third identity
in \eqref{eq:radial-equality} gives
$Y=f_\nu\nu=\nu/H_\phi$. Hence
\eqref{eq:radial-position} yields
\[
x=\frac{V}{H_\phi}\nu\qquad\text{on }\Sigma.
\]
For every tangent vector $X\in T\Sigma$,
\[
X(|x|^2)=2\langle x,X\rangle
=2\frac{V}{H_\phi}\langle\nu,X\rangle=0.
\]
Thus $|x|$ is constant on each connected boundary
component, which is therefore a sphere centered at
the origin. Since $V/H_\phi>0$, the outward normal
is $\nu=x/|x|$, excluding inner boundary components.
As $\Omega$ is bounded and connected, we conclude that
$\Omega=B_R(0)$ for some $R>0$.

Furthermore, \eqref{eq:radial-position} and the second
identity in \eqref{eq:radial-equality} give
\[
\overline\nabla\phi
=-\frac{|\overline\nabla\phi|^2}{V-m}\,x.
\]
Thus $\phi$ is constant on every sphere centered at
the origin, and therefore radial throughout $\Omega$.
By \eqref{NV}, $V$ is radial as well.

Conversely, suppose $\Omega=B_R(0)$. Since
$x=R\nu$ on $\Sigma$, \eqref{NV}
gives $V=RH_\phi$ there. Moreover,
\[
\overline\diver_\phi x
=m+1-\langle x,\overline\nabla\phi\rangle
=V+1=\mathcal N,
\qquad
\frac{\mathcal N}{\mathcal N-1}V=\mathcal N.
\]
The weighted divergence theorem yields
\begin{align*}
\int_\Sigma\frac{V}{H_\phi}\,d\sigma_\phi
&=R\int_\Sigma d\sigma_\phi
=\int_\Sigma\langle x,\nu\rangle\,d\sigma_\phi\\
&=\int_\Omega\mathcal N\,d\mu_\phi
=\int_\Omega\frac{\mathcal N}{\mathcal N-1}
V\,d\mu_\phi.
\end{align*}
Thus equality holds in \eqref{eq:general-HK}.
\end{proof}

\bibliographystyle{amsplain}
\bibliography{references}
\end{document}